\documentclass{amsart}%
\usepackage{amssymb}
\usepackage{amsmath}
\usepackage{graphicx}
\usepackage{amsfonts}
\usepackage{color}
\newtheorem{theorem}{Theorem}
\theoremstyle{plain}

\newtheorem{definition}{Definition}
\newtheorem{example}{Example}

\newtheorem{lemma}{Lemma}

\newtheorem{remark}{Remark}

\numberwithin{equation}{section}

\input{tcilatex}

\begin{document}
\title{Some Naturally Occurring Examples of $A_{\infty}$-bialgebras}
\author{Ainhoa Berciano$^{1}$}
\address{Departamento de Did\'actica de la Matem\'atica y de las Ciencias Experimentales\\
Universidad del Pa\'{\i}s Vasco\\
Bilbao, C.P. 48014}
\author{Ronald Umble$^{2}$}
\address{Department of Mathematics\\
Millersville University of Pennsylvania\\
Millersville, PA. 17551}
\thanks{$^{1}$Work partially supported by project \textquotedblleft EHU09/04\textquotedblright and by the PAICYT research project FQM-296}
\thanks{$^{2}$This research was funded in part by a Millersville University faculty
research grant}
\date{July 8, 2010}

\begin{abstract}
Let $p$ be an odd prime. We show that when $n\geq 3$, each tensor factor
$E\otimes\Gamma$ of $H_{\ast}\left(  \mathbb{Z},n;\mathbb{Z}_{p}\right)  $ is
an $A_{\infty}$-bialgebra with non-trivial structure. We give explicit
formulas for the structure maps and the relations among them. Thus
$E\otimes\Gamma$ is a naturally occurring $A_{\infty}$-bialgebra.

\end{abstract}
\maketitle

\section{Introduction}

Let $p$ be an odd prime and let $n \geq3$. Let $E\left(  v,2n+1\right)  $
denote the $\mathbb{Z}_{p}$-exterior algebra on a generator $v$ of dimension
$2n+1$ and let $\Gamma(w,2np+2)$ denote the $\mathbb{Z}_{p}$-divided power
algebra on a generator $w$ of dimension $2np+2$. In \cite{Cartan} and
\cite{EM54}, Cartan, Eilenberg and Mac Lane showed that $H_{\ast}\left(
\mathbb{Z},n;\mathbb{Z}_{p}\right)  $ factors as an infinite tensor product
with infinitely many tensor factors of the form $E \otimes\Gamma$; for
example,
\[
H_{\ast}\left(  \mathbb{Z},3;\mathbb{Z}_{p}\right)  \approx\bigotimes_{i\geq
0}E(v_{i},2p^{i}+1)\otimes\Gamma(w_{i},2p^{i+1}+2).
\]

The main result in this paper is that for each $i\geq0$, the factor $$A=E\left(
v_{i},2m+1\right) \linebreak \otimes\Gamma(w_{i},2mp+2)\subset H_{\ast}\left(
\mathbb{Z},n;\mathbb{Z}_{p}\right)  $$ is an $A_{\infty}$-bialgebra with
exactly three non-trivial \textquotedblleft structurally
compatible\textquotedblright operations, namely, a multiplication $\mu$, a
comultiplication $\Delta$, and an operation $\Delta_{p}:A\rightarrow
A^{\otimes p}$ of degree $p-2$. The operations $\Delta$ and $\Delta_{p}$
define the $A_{\infty}$-coalgebra structure of $A$ obtained by the first
author using techniques of homological perturbation theory (see \cite{BR07}),
and can be realized as a contraction of the reduced bar construction $B\left(
\mathbb{Z}\left[  u\right]  /\left(  u^{p}\right)  \right)  $ with $\left\vert
u\right\vert =2n.$

Since $A$ is a Hopf algebra, $\Delta$ and $\mu$ are compatible in the sense
that $\Delta$ is an algebra map, and it is natural to ask whether $\Delta_{p}$
and $\mu$ are compatible in some analogous way. Indeed, $$f^{p}=\left(
\Delta\otimes\mathbf{1}^{\otimes p-2}\right)  \cdots\left(  \Delta
\otimes\mathbf{1}\right)  \Delta$$ is an algebra map, and $\Delta_{p}$ is
compatible with $\mu$ as an \textquotedblleft$\left(  f^{p},f^{p}\right)
$-derivation of degree $p-2.$\textquotedblright\ Consequently, $A$ is an
$A_{\infty}$-bialgebra as defined by S. Saneblidze and the second author in
\cite{SU3}. We will refer to $A$ as a \textquotedblleft Hopf $A_{\infty}%
$-coalgebra\textquotedblright\ to emphasize the compatibility of $\Delta$ and
$\Delta_{p}$ with $\mu$.

The paper is organized as follows: Section 2 reviews the notion
of an $A_{\infty}$-(co)algebra and the related tilde (co)bar
construction. Section 3 reviews the $A_{\infty}$-coalgebra
structure on $E\otimes\Gamma$ mentioned above. Section 4 reviews
the definition of the S-U diagonals $\Delta_{P}$ and $\Delta_{K}$ on the cellular
chains of permutahedra and associahedra given in
\cite{SU2}. In Section 5 we give an exposition of the general
notion of a \textquotedblleft higher derivation\textquotedblright,
i.e., a $\Delta$-derivation with respect to a $\Delta$-compatible
family of maps indexed by the faces of a family of polytopes
$X=\sqcup_{n\geq0}X_{n}$; the ideas in this section are implicit
in \cite{SU3}.  We conclude the paper with Section 6, in which we
prove Theorem 4: \textit{Let} $p$ \textit{be an odd prime, let} $i \geq 0$, \textit{and let} $n\geq 1$. \textit{Then}
$$E(v_{i},2np^i+1)\otimes\Gamma(w_{i},2np^{i+1}+2)$$ \textit{is a Hopf} $A_{\infty}%
$\textit{-coalgebra over} $\mathbb{Z}_{p}.$

\section{$A_{\infty}$-Coalgebras and Related Constructions\label{Section 2}}

Let $R$ be a commutative ring with unity and let $M$ be a graded $R$-module.
Let $\overline{M}=M/M_{0},$ and let $\uparrow$ and $\downarrow$ denote the
suspension and desuspension operators, which shift dimensions $+1$ and $-1,$
respectively. Given a map $f:M\rightarrow N$ of graded modules, let $f_{i,j}$
denote the map $$\mathbf{1}^{\otimes i}\otimes f\otimes\mathbf{1}^{\otimes
j}:N^{\otimes i}\otimes M\otimes N^{\otimes j}\rightarrow N^{\otimes i+j+1}.$$

Given a connected differential graded (DG)\ algebra $\left(  A,d,\mu\right)
$, the \emph{bar construction of} $A$ is the tensor coalgebra $BA=T^{c}\left(
\uparrow\bar{A}\right)  $ with cofree coproduct
\[
\Delta_{B}[a_{1}|\cdots|a_{n}]=1\otimes\lbrack a_{1}|\cdots|a_{n}%
]+[a_{1}|\cdots|a_{n}]\otimes1+\sum_{i=1}^{n-1}[a_{1}|\cdots|a_{i}%
]\otimes\lbrack a_{i+1}|\cdots|a_{n}]
\]
and differential $d_{B}=d_{t}+d_{a}$, where%
\[
d_{t}=\sum_{i=1}^{n}\left(  \uparrow d\downarrow\right)  _{i-1,n-i}\text{
\ and \ }d_{a}=\sum_{i=1}^{n-1}\left(  \uparrow\mu\downarrow^{\otimes
2}\right)  _{i-1,n-i-1}.
\]
Dually, given a simply connected DG coalgebra $\left(  C,d,\Delta\right)  $,
the \emph{cobar construction of} $C$ is the tensor algebra $\Omega
C=T^{a}\left(  \downarrow\bar{C}\right)  $ with differential $d_{\Omega}%
=d_{t}+d_{c},$ where
\[
d_{c}=\sum_{i=1}^{n}\left(  \downarrow^{{\otimes}2}\Delta\uparrow\right)
_{i-1,n-i}.
\]

The notion of an $A_{\infty}$\emph{-algebra} was defined by J. Stasheff in his
seminal paper \cite{Sta63}; here we review the dual notion of an $A_{\infty}%
$\emph{-coalgebra}. Given a simply connected graded $R$-module $A$ together
with a family of $R$-multilinear operations $\{\psi^{k}\in Hom^{k-2}\left(
A,A^{\otimes k}\right)  \}_{k\geq1},$ let
\[
d=\sum\limits_{\substack{0\leq i\leq n-1;\text{ }1\leq k\leq n\\n\geq
1}}\left(  \downarrow^{\otimes k}\psi^{k}\uparrow\right)  _{i,n-i-1}:\Omega
A\rightarrow\Omega A.
\]
Then $\left(  A,\psi^{n}\right)  _{n\geq1}$ is an $A_{\infty}$%
-\emph{coalgebra} if $d\circ d=0,$ and the structure
relations%
\begin{equation}
\sum_{\substack{0\leq i\leq n-j-1\\0\leq j\leq n-1;\text{ }n\geq1}}\left(
-1\right)  ^{j\left(  n+i+1\right)  }\psi_{i,n-i-j-1}^{j+1}\psi^{n-j}=0
\label{relations}%
\end{equation}
are given by the homogeneous components of the equation $d\circ d=0.$ The
signs in (\ref{relations}) are the Koszul signs that appear when factoring out
suspension and desuspension operators. Once the signs have been determined, we
drop the simple-connectivity assumption and obtain the following abstract
definition: An $A_{\infty}$\emph{-coalgebra }is a graded $R$-module\ $A$
together with a family of multilinear operations $$\{\psi^{k}\in Hom^{k-2}%
\left(  A,A^{\otimes k}\right)  \}_{k\geq1}$$ that satisfy the structure
relations in (\ref{relations}).

\section{The $A_{\infty}$-coalgebra $E\otimes\Gamma$}

For $n\in\mathbb{N}$, let $Q_{p}\left(  u,2n\right)  =\mathbb{Z}\left[
u\right]  /\left(  u^{p}\right)  \ $with $\left\vert u\right\vert =2n.$ Recall
that the exterior algebra $E(v,2n-1)$ on a generator $v$ of degree $2n-1$ is a
Hopf algebra with primi\-tively generated coproduct $\Delta$. As a module, the
divided power algebra $\Gamma(w,2n)$ is generated by ${\gamma_{i}=\gamma
_{i}(w)}_{i\geq1},$ where $\gamma_{1}(w)=w$, and the algebra structure is
defined by
\[
{\gamma_{i}\gamma_{j}=\frac{(i+j)!}{i!j!}\gamma_{i+j}}.
\]
Furthermore, $\Gamma(w,2n)$ is a Hopf algebra with respect to the coproduct
generated by $$\Delta(\gamma_{k}(u))=\sum_{i+j=k}\gamma_{i}(u)\otimes\gamma
_{j}(u).$$

In \cite{Pro84}, A. Prout\'{e} took the first steps toward computing the
$A_{\infty}$-coalgebra structure on $H=H_{\ast}\left(  K(\pi,n);\mathbb{Z}%
_{p}\right)  $. He showed that $H$ is a classical coalgebra when $p=2$. But
when $p$ is an odd prime, the $A_{\infty}$-coalgebra structure explodes and he
only obtained partial results in certain special cases. Thanks to Eilenberg
and Mac Lane \cite{EM54}, there is a contraction (a special type of chain
homotopy equivalence)
\[
B(Q_{p}(u_{i},2np^{i}))\rightarrow E(v_{i},2np^{i}+1)\otimes\Gamma
(w_{i},2np^{i+1}+2)
\]
by which we obtain the decomposition
\[
H_{\ast}(\mathbb{Z},3,\mathbb{Z}_{p})\approx\bigotimes_{i\geq0}E\left(
w_{i},2p^{i}+1\right)  \otimes\Gamma\left(  x_{i},2p^{i+1}+2\right)  
\]
(see Theorem 4.5 in \cite{BR07}). This idea extends inductively to $H_{\ast
}(\mathbb{Z},n,\mathbb{Z}_{p})$ (see Theorem 4.24 of \cite{BR07}). When $n=4$,
let $n_{i}=p^{i+1}+1$; then
\[
H_{\ast}\left(  \mathbb{Z},4;\mathbb{Z}_{p}\right)  \approx\bigotimes_{i\geq0}\left[
\Gamma\left(  y_{i},2p^{i}+2\right)  \otimes\bigotimes_{j\geq0}E\left(
a_{j},2n_{i}p^{j}+1\right)  \otimes\Gamma\left(  b_{j},2n_{i}p^{j+1}+2\right)
\right]  ;
\]
when $n=5$ we have
\[
H_{\ast}\left(  \mathbb{Z},5;\mathbb{Z}_{p}\right)  \approx\bigotimes_{i\geq0}\left[
\left(  \bigotimes_{k\geq0}E_{i,k}\otimes\Gamma_{i,k}\right)  \otimes
\bigotimes_{j\geq0}\Gamma_{i,j}\otimes\left(  \bigotimes_{l\geq0}%
E_{i,j,l}\otimes\Gamma_{i,j,l}\right)  \right]  ;
\]
and so on. In \cite{BR07}, the first author used these decompositions to
extend Prout\'{e}'s results and obtain

\begin{theorem}
For $n\geq3,$ the induced $A_{\infty}$-coalgebra operation $\Delta_{q}$ on
$H_{\ast}(K(\mathbb{Z},n);\mathbb{Z}_{p})$ vanishes whenever $q\neq i(p-2)+2$
and $i\geq0$.
\end{theorem}

This result follows immediately from

\begin{theorem}
\label{berciano}For all $m\in\mathbb{N}$ and every odd prime $p,$ the
$\mathbb{Z}_{p}$-Hopf algebra $$A=E(v,2m+1)\otimes\Gamma(w,2mp+2)$$ is a
non-trivial $A_{\infty}$-coalgebra. The induced operation $\Delta
_{q}:A\rightarrow A^{\otimes q}$ is non-trivial if and only if $q=2$ or $q=p.$
In fact, for $i=0,1$ and $\gamma_{j}=\gamma_{j}\left(  w\right)  $ we have%
\begin{equation}
\Delta_{2}(v^{i}\gamma_{j})=\sum_{k=0}^{i}\sum_{l=0}^{j}v^{k}\gamma_{l}\otimes
v^{i-k}\gamma_{j-l}; \label{Delta2}%
\end{equation}%
\begin{equation}
\Delta_{p}(v^{i}\gamma_{j})=\sum_{k_{1}+\cdots+k_{p}=j-1}v^{i+1}\gamma_{k_{1}%
}\otimes\cdots\otimes v^{i+1}\gamma_{k_{p}}. \label{Deltap}%
\end{equation}

\end{theorem}

The coproduct $\Delta_{2}$ defined in (\ref{Delta2}) is the induced coproduct
on the tensor product of coalgebras and is compatible with the induced
multiplication $\mu$ as an algebra map. And as we shall see, $\Delta_{p}$ is
also compatible with $\mu$ as an \textquotedblleft$\left(  f^{p},f^{p}\right)
$-derivation of degree $p.$\textquotedblright

\section{The S-U Diagonals on Permutahedra and Associahedra}

This section gives a brief review of the S-U diagonals $\Delta_{P}$  and $\Delta_{K}$ on
permutahedra $P=\sqcup_{n \geq 1} P_n$ and associahedra $K=\sqcup_{n \geq 2} K_{n} $ (up to sign); for details see \cite{SU2}.
Alternative constructions of $\Delta_{K}$ were subsequently given by Markl and
Shnider \cite{MS} and Loday \cite{Loday}.

Let $n\in\mathbb{N}$ and let $\underline{n}=\{1,2,\dots,n\}.$ A matrix $E$
with entries from $\left\{  0\right\}  \cup\underline{n}$ is a \emph{step
matrix} if the following conditions hold:

\begin{itemize}
\item Each element of $\underline{n}$ appears as an entry of $E$ exactly once.

\item Elements of $\underline{n}$ in each row and column of $E$ form an
increasing contiguous block.

\item Each diagonal parallel to the main diagonal of $E$ contains exactly one
element of $\underline{n}.$
\end{itemize}

\noindent The non-zero entries in a step matrix form a continuous staircase
connecting the lower-left and upper-right entries. There is a bijective
correspondence between step matrices with non-zero entries in $\underline{n}$
and permutations of $\underline{n}$.

Given a $q\times p$ integer matrix $M=\left(  m_{ij}\right)  ,$ choose proper
subsets $S_{i}\subset\{  \text{non-zero} \linebreak \text{entries in }\operatorname*{row}%
\left(  i\right)  \}  $ and $T_{j}\subset\left\{  \text{non-zero entries
in }\operatorname*{col}\left(  j\right)  \right\}  ,$ and define
\emph{down-shift} and \emph{right-shift} operations\emph{\ }$D_{S_{i}}$ and
$R_{T_{j}}$ on $M$ as follows:

\begin{enumerate}
\item[\textit{(i)}] If $S_{i}\neq\varnothing,$ $\max\operatorname*{row}%
(i+1)<\min S_{i}=m_{ij},$ and $m_{i+1,k}=0$ for all $k\geq j$, then $D_{S_{i}%
}M$ is the matrix obtained from $M$ by interchanging each $m_{ik}\in S_{i}$
with $m_{i+1,k};$ otherwise $D_{S_{i}}M=M.$

\item[\textit{(ii)}] If $T_{j}\neq\varnothing,$ $\max\operatorname*{col}%
(j+1)<\min T_{j}=m_{ij},$ and $m_{k,j+1}=0$ for all $k\geq i,$ then\emph{\ }%
$R_{T_{j}}M$ is the matrix obtained from $M$ by interchanging each $m_{k,j}\in
T_{j}$ with $m_{k,j+1};$ otherwise $R_{Tj}M=M.$
\end{enumerate}

\noindent Given a $q\times p$ step matrix $E\ $together with subsets
$S_{1},\ldots,S_{q}$ and $T_{1},\ldots,T_{p}$ as above, there is the
\emph{derived matrix}%
\[
R_{T_{p}}\cdots R_{T_{2}}R_{T_{1}}D_{S_{q}}\cdots D_{S_{2}}D_{S_{1}}E.
\]
In particular, step matrices are derived matrices under the trivial action
with $S_{i}=T_{j}=\varnothing$ for all $i,j$.

Let $a=A_{1}|A_{2}|\cdots|A_{p}$ and $b=B_{q}|B_{q-1}|\cdots|B_{1}$ be
partitions of $\underline{n}.$ Then $a\times b$ is a $(p,q)$%
\emph{-complementary pair} (CP) if there is a $q\times p$ derived matrix
$M=\left(  m_{ij}\right)  $ such that $A_{j}=\{m_{ij}\neq0\mid1\leq i\leq q\}$
and $B_{i}=\{m_{ij}\neq0\mid1\leq j\leq p\}.$ Thus $\left(  p,q\right)  $-CPs,
which are in one-to-one correspondence with derived matrices, identify a
particular set of product cells in $P_{n}\times P_{n}$.

\begin{definition}
Define $$\Delta_{P}:C_{0}\left(  P_{1}\right)  \rightarrow C_{0}\left(
P_{1}\right)  \otimes C_{0}\left(  P_{1}\right)  $$ by $\Delta_{P}%
(\underline{1})=\underline{1}\otimes\underline{1}$. Inductively, having
defined $$\Delta_{P}:C_{\ast}\left(  P_{k}\right)  \rightarrow C_{\ast}\left(
P_{k}\right)  \otimes C_{\ast}\left(  P_{k}\right)  $$ for all $k\leq n$,
define $\Delta_{P}$ on $\underline{n+1}\in C_{n}(P_{n+1})$ by
\[
\Delta_{P}(\underline{n+1})=\sum_{\substack{(p,q)\text{-CPs }a\times
b\\p+q=n+2}}\pm\text{ }a\otimes b
\]
and extend $\Delta_{P}$ to all of $C_{\ast}\left(  P_{n+1}\right)  $
multiplicatively, i.e., define $\Delta_{P}$ on a generator $u_{1}|\cdots
|u_{r}\in C_{n-r+1}\left(  P_{n+1}\right)  $ by $$\Delta_{P}\left(
u_{1}|\cdots|u_{r}\right)  =\Delta_{P}\left(  u_{1}\right)  |\cdots|\Delta
_{P}\left(  u_{r}\right)  .$$
\end{definition}

Recall that faces of $P_{n}$ in codimension $k$ are indexed by
planar-rooted trees with $n+1$ leaves and $k+1$ levels (PLTs).

\begin{example}
In terms of PLTs, the diagonal $\Delta_{P}$ on the top cell of $P_{3}$ (up to sign) is given by
\begin{equation*}
\FRAME{itbpF}{2.5668in}{0.9314in}{0in}{}{}{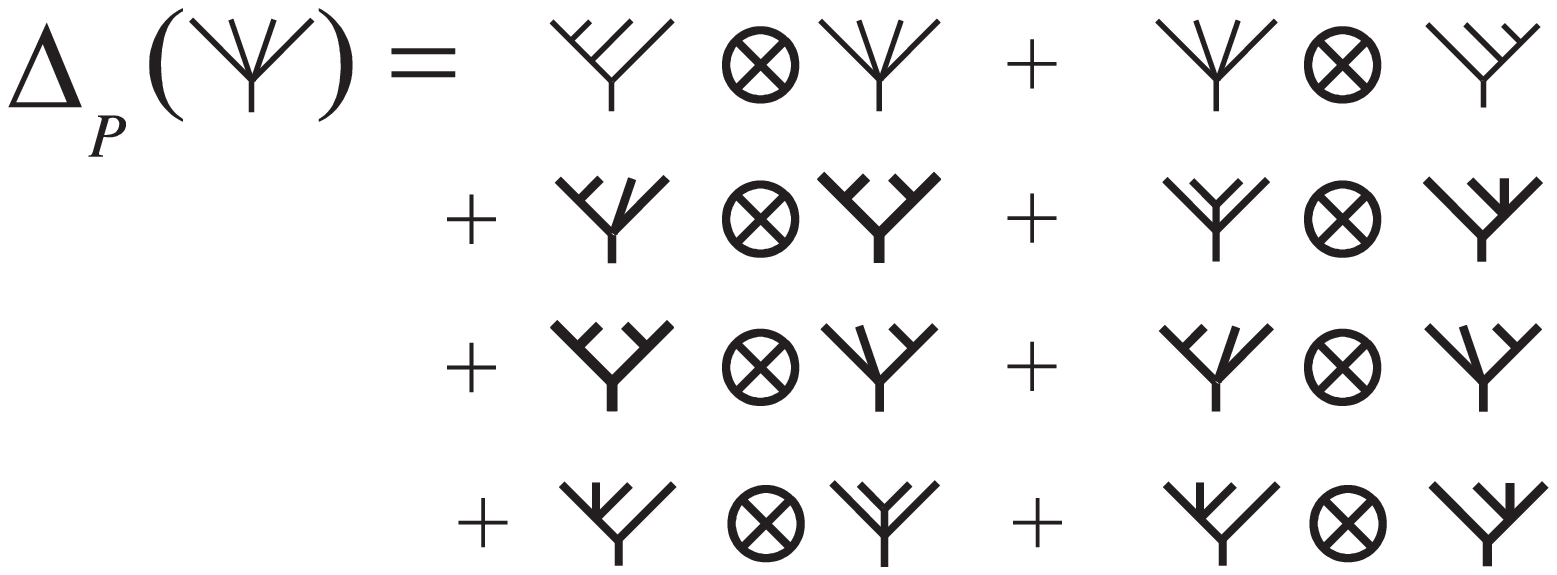}{\raisebox{-0.9314in}{\includegraphics[height=0.9314in]{deltap3.eps}}}.
\end{equation*}
\end{example}

The diagonal $\Delta_P$ descends to a diagonal $\Delta_K$ on $C_*(K)$ via Tonks' cellular projection $\vartheta_0 : P_n \to K_{n+1}$, which forgets levels (see \cite{Tonks}). Faces of $P_{n}$ indexed by PLTs
with multiple nodes in the same level degenerate under $\vartheta_{0}$, and
corresponding generators span the kernel of the induced map $\vartheta
_{0}:C_{\ast}\left(  P_{n}\right)  \rightarrow C_{\ast}\left(  K_{n+1}\right)
$. The diagonal $\Delta_{K}$ is given by
\[
\Delta_{K}\vartheta_{0}=(\vartheta_{0}\otimes\vartheta_{0})\Delta_{P}.
\]

\begin{example}
When $n=3,$ the components $1|23\otimes 13|2$ and $13|2\otimes 3|12$ of $%
\Delta _{P}\left( \underline{3} \right) $ degenerate under $\vartheta$ because the tree
corresponding to $13|2$ has two vertices in level $1$; equivalently, $%
\dim \left( 13|2\right) =1$ whereas $\dim \vartheta \left( 13|2\right) =0.$
Therefore (up to sign) the diagonal on the top cell of $K_{4}$ is given by
\begin{equation*}
\FRAME{itbpF}{2.45in}{0.70in}{0in}{}{}{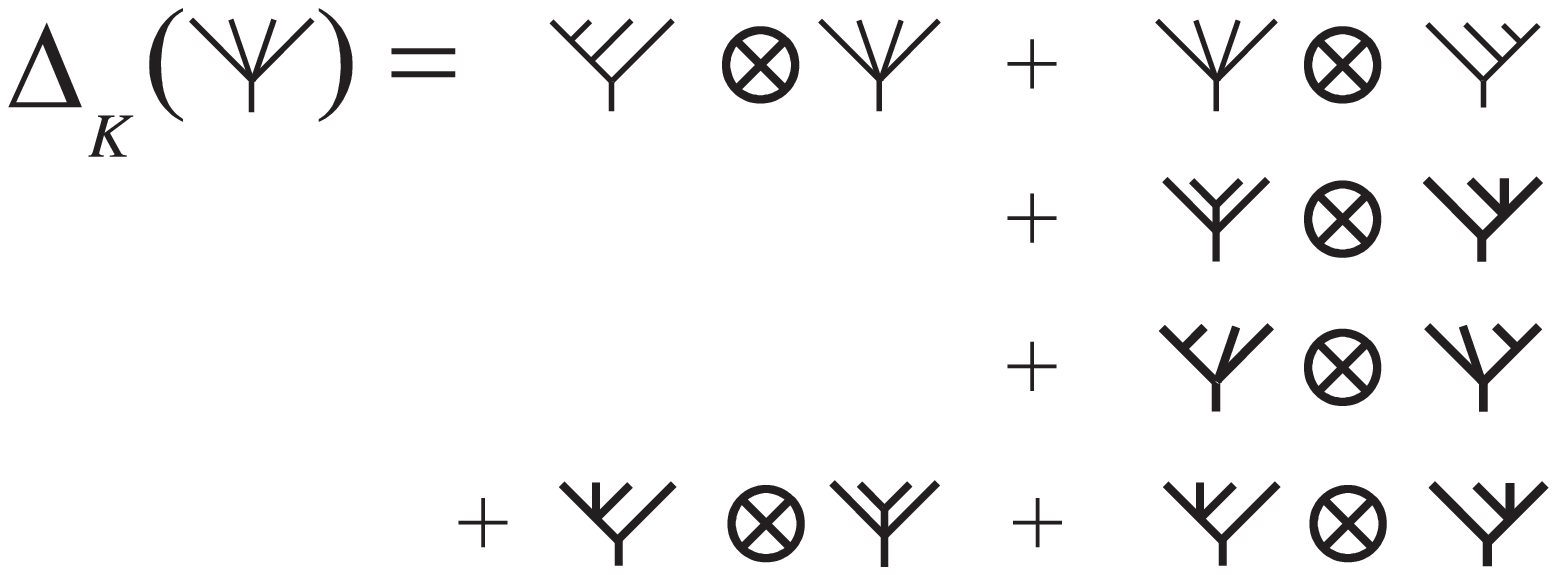}{\raisebox{-0.70in}{\includegraphics[height=0.70in]{deltak4.eps}}} \hspace*{3.0mm}.
\end{equation*}
\end{example}

\section{$\Delta$-derivations and $\Delta_{X}$-compatible families}

Let $\left\{X_{n}\right\}_{n\geq 0}$ be a family of contractible polytopes such that $%
\dim X_{n}=n,$ let $X=\sqcup _{n\geq 0}X_{n}$ and assume that the cellular
chains $C_{\ast }\left( X\right) $ are equipped with a diagonal
approximation $$\Delta _{X}:C_{\ast }\left( X\right) \rightarrow C_{\ast
}\left( X\right) \otimes C_{\ast }\left( X\right) .$$ In this section we
introduce the general notion of a $\Delta $-derivation homotopy with respect to a $
\Delta _{X}$-compatible family of maps. When $X_{n}$ is the $n$-simplex $s_{n}$,
our definition agrees with the notion of a \emph{high derivation} defined by
T. Kadeishvili in \cite{Kadeishvili2}. When $X_{n}$ is the $n$-dimensional
permutahedron $P_{n+1}$, the notion of a $\Delta $%
-derivation with respect to a $\Delta _{P}$-compatible family is encoded in
the construction of the biderivative given by S. Saneblidze and the second
author (see \cite{SU3}, \cite{SU4} and \cite{Umble}).

For each $n,$ let $n_{k}$ be the number of $k$-faces of $X_{n}$ and choose a
system of generators $\left\{x_{i}^{k}\right\}_{0\leq k\leq n;1\leq i\leq
n_{k}}$ for $C_{\ast}\left( X_{n}\right) .$ Let $X_{i}^{k}$ denote the
smallest subcomplex of $X_{n}$ containing the $k$-face associated with $%
x_{i}^{k}.$ Given DGAs $\left( A,\mu_{A},d_{A}\right) $ and $\left(
B,\mu_{B},d_{B}\right) $, let $\Theta:C_{\ast}\left( X_{n}\right)
\rightarrow Hom\left( A,B\right) $ be a map of degree zero and let $%
\Theta_{i}^{k}=\Theta|_{C_{\ast}\left( X_{i}^{k}\right) }$.

\begin{definition}
\textit{The family of maps}
\begin{equation*}
\mathfrak{F}_{n}=\left\{ \Theta \left( x_{i}^{k}\right) \right\}_{0\leq
k<n;1\leq i\leq n_{k}}
\end{equation*}
\textit{is } \textbf{$\Delta _{X}$-compatible}\textit{\ if each }$%
\Theta _{i}^{k}$ \textit{is a chain map and the following diagram commutes}:
\begin{equation*}
\begin{array}{ccc}
C_{k}(X_{i}^{k}) & \overset{\Delta _{X}}{\longrightarrow } &
\sum\limits_{p+q=k}C_{p}(X_{i}^{k})\otimes C_{q}(X_{i}^{k}) \\
&  &  \\
\text{\ \ }_{{\LARGE \Theta }_{i}^{k}}\text{ }\downarrow &  & \text{\ \ \ \
\ \ \ \ \ }\downarrow \text{ }_{{\LARGE \Theta }_{i}^{k}{\LARGE \otimes }%
\text{ }{\LARGE \Theta }_{i}^{k_{\text{ \ }}}} \\
&  &  \\
Hom^{k}\left( A,B\right) &  & \sum\limits_{p+q=k}Hom^{p}\left( A,B\right)
\otimes Hom^{q}\left( A,B\right) \\
&  &  \\
_{\left( {\LARGE \mu }_{A}\right) ^{\ast }}\text{ }\downarrow &  &
\downarrow \text{ \ }\approx \\
&  &  \\
Hom^{k}\left( A\otimes A,B\right) & \underset{\left( {\LARGE \mu }%
_{B}\right) _{\ast }}{\longleftarrow } & Hom^{k}\left( A\otimes A,B\otimes
B\right) .%
\end{array}%
\end{equation*}%
Let $\mathfrak{F}_{n}$\ be a $\Delta _{X}$\textit{-compatible family of
maps. }The map $T=\Theta \left( x_{1}^{n}\right) :A\rightarrow B$ associated
with the top dimensional cell of $X$ is a \textbf{$\Delta $-derivation
with respect to $\mathfrak{F}_{n}$}\ if the diagram above commutes when $k=n$%
. If in addition, $\Theta $ is a chain map, then $T$ is a \textbf{$\Delta
$-derivation homotopy with respect to} $\mathfrak{F}_{n}$. There is the dual
notion of a $\Delta $-coderivation homotopy with respect to a $\Delta _{X}$%
-compatible family.
\end{definition}

Let $\sigma_{n,2}:\left(A^{\otimes n}\right)^{\otimes 2} \to \left(A^{\otimes 2}\right)^{\otimes n}$ be the canonical permutation of tensor factors.  Then for example,
$$ \sigma_{3,2}\left(a_1|a_2|a_3\otimes b_1|b_2|b_3\right)=\left(-1\right)^{|b_1||a_2|+|b_1||a_3|+|b_2||a_3|}a_1|b_1 \otimes a_2|b_2 \otimes a_3|b_3.$$

\begin{example}
Set $X_{n}=K_{n+2}$ and let $\Delta _{K}$ be the S-U diagonal on
$C_*\left( K \right)$. Given a DGA $\left( A,d,\mu \right) ,$  let $d^{\otimes i}$ denote the free
linear extension of $d$ to $A^{\otimes i}$.  Choose an arbitrary family
of DG module maps $\{\Delta _{i}\in Hom^{i-2}\left(
A,A^{\otimes i}\right) \}_{i\geq 2}$. For notational simplicity, identify
$\Delta_{i}$ with the down-rooted $i$-leaf corolla and other
down-rooted planar trees with their
corresponding compositions in $Hom\left( A,A^{\otimes \ast }\right) $. When $n=2,$ $K_{2}$ is a point. If \
\FRAME{itbpF}{0.1358in}{0.1237in}{0.0138in}{}{}{t2.eps}{\raisebox{-0.1237in}{\includegraphics[height=0.1237in]{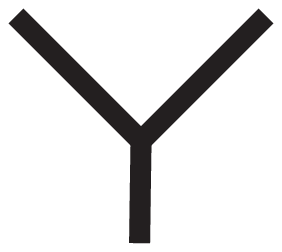}}} \
is $\Delta $-derivation homotopy with respect to the empty family $\mathfrak{F}_{2},$ then \
\FRAME{itbpF}{0.1358in}{0.1237in}{0.0208in}{}{}{t2.eps}{\raisebox{-0.1237in}{\includegraphics[height=0.1237in]{T2.eps}}} is a DGA map, i.e.,
\begin{equation*}
d^{\otimes 2}\text{}
\FRAME{itbpF}{0.1358in}{0.1237in}{0.0138in}{}{}{t2.eps}{\raisebox{-0.1237in}{\includegraphics[height=0.1237in]{T2.eps}}}-
\FRAME{itbpF}{0.1358in}{0.1237in}{0.0138in}{}{}{t2.eps}{\raisebox{-0.1237in}{\includegraphics[height=0.1237in]{T2.eps}}}
\text{ }d=0\text{ \ and \ }
\FRAME{itbpF}{0.1358in}{0.1237in}{0.0208in}{}{}{t2.eps}{\raisebox{-0.1237in}{\includegraphics[height=0.1237in]{T2.eps}}}
\text{ }\mu =\mu \text{ }\sigma_{2,2}\left( \text{ }
\FRAME{itbpF}{0.1358in}{0.1237in}{0.0208in}{}{}{t2.eps}{\raisebox{-0.1237in}{\includegraphics[height=0.1237in]{T2.eps}}}
\text{ }\otimes \text{ }
\FRAME{itbpF}{0.1358in}{0.1237in}{0.0138in}{}{}{t2.eps}{\raisebox{-0.1237in}{\includegraphics[height=0.1237in]{T2.eps}}}
\text{ }\right) .
\end{equation*}%
When $n=3,$ $K_{3}$ is an interval. If $\mathfrak{F}_{3}=
\left\{
\text{ }
\FRAME{itbpF}{0.1436in}{0.1418in}{0.0208in}{}{}{t21-2.eps}{\raisebox{-0.1418in}{\includegraphics[height=0.1418in]{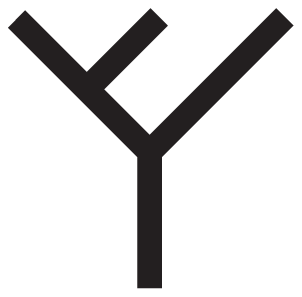}}},
\text{ } \FRAME{itbpF}{0.1436in}{0.1418in}{0.0277in}{}{}{t12-2.eps}{\raisebox{-0.1418in}{\includegraphics[height=0.1418in]{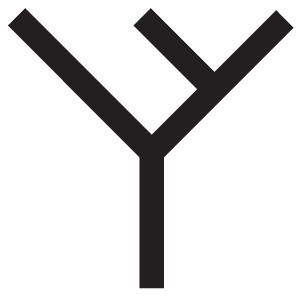}}}
\right\}
$ is a $\Delta _{K}$-compatible family of compositions, then \
\FRAME{itbpF}{0.1436in}{0.1418in}{0.0208in}{}{}{t21-2.eps}{\raisebox{-0.1418in}{\includegraphics[height=0.1418in]{T21-2.eps}}} and
\FRAME{itbpF}{0.1436in}{0.1418in}{0.0277in}{}{}{t12-2.eps}{\raisebox{-0.1418in}{\includegraphics[height=0.1418in]{T12-2.eps}}}
are DGA\ maps. If
\FRAME{itbpF}{0.1574in}{0.1427in}{0.0277in}{}{}{t3.eps}{\raisebox{-0.1427in}{\includegraphics[height=0.1427in]{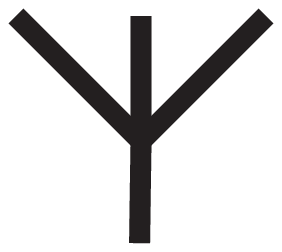}}}
is a $\Delta $-derivation homotopy with respect to $\mathfrak{F}_{3},$ then
\FRAME{itbpF}{0.1574in}{0.1427in}{0.0138in}{}{}{t3.eps}{\raisebox{-0.1427in}{\includegraphics[height=0.1427in]{T3.eps}}}
\ is a $\left( \text{ }
\FRAME{itbpF}{0.1436in}{0.1418in}{0.0277in}{}{}{t21-2.eps}{\raisebox{-0.1418in}{\includegraphics[height=0.1418in]{T21-2.eps}}}
\text{ \ },\text{ }
\FRAME{itbpF}{0.1436in}{0.1418in}{0.0277in}{}{}{t12-2.eps}{\raisebox{-0.1418in}{\includegraphics[height=0.1418in]{T12-2.eps}}}\text{ }\right) $-derivation homotopy, i.e.,
\begin{equation*}
d^{\otimes 3}\text{ }
\FRAME{itbpF}{0.1574in}{0.1427in}{0.0277in}{}{}{t3.eps}{\raisebox{-0.1427in}{\includegraphics[height=0.1427in]{T3.eps}}}
\text{ }+\text{ }
\FRAME{itbpF}{0.1574in}{0.1427in}{0.0277in}{}{}{t3.eps}{\raisebox{-0.1427in}{\includegraphics[height=0.1427in]{T3.eps}}}
\text{ }d=\text{ }
\FRAME{itbpF}{0.1436in}{0.1418in}{0.0415in}{}{}{t12-2.eps}{\raisebox{-0.1418in}{\includegraphics[height=0.1418in]{T12-2.eps}}}-
\FRAME{itbpF}{0.1436in}{0.1418in}{0.0346in}{}{}{t21-2.eps}{\raisebox{-0.1418in}{\includegraphics[height=0.1418in]{T21-2.eps}}}
\text{ \ and \ }
\end{equation*}
\begin{equation*}
\FRAME{itbpF}{0.1574in}{0.1427in}{0.0277in}{}{}{t3.eps}{\raisebox{-0.1427in}{\includegraphics[height=0.1427in]{T3.eps}}}
\mu =\mu ^{\otimes 3}\text{ }\sigma _{3,2}\left( \text{ }
\FRAME{itbpF}{0.1436in}{0.1418in}{0.0415in}{}{}{t21-2.eps}{\raisebox{-0.1418in}{\includegraphics[height=0.1418in]{T21-2.eps}}}
\text{ }\otimes \text{ }
\FRAME{itbpF}{0.1574in}{0.1427in}{0.0277in}{}{}{t3.eps}{\raisebox{-0.1427in}{\includegraphics[height=0.1427in]{T3.eps}}}
\text{ }+\text{ }
\FRAME{itbpF}{0.1574in}{0.1427in}{0.0277in}{}{}{t3.eps}{\raisebox{-0.1427in}{\includegraphics[height=0.1427in]{T3.eps}}}
\text{ }\otimes \text{ }\FRAME{itbpF}{0.1436in}{0.1418in}{0.0346in}{}{}{t12-2.eps}{\raisebox{-0.1418in}{\includegraphics[height=0.1418in]{T12-2.eps}}}\text{ }\right) .
\end{equation*}
In the case of the pentagon $K_{4},$ assume that $\mathfrak{F}_{4}=\left\{
\Theta \left( x_{i}^{k}\right) \right\} _{k=0,1}$ is a $\Delta _{K}$%
-compatible family. Then
\begin{equation*}
\begin{tabular}{l}
\FRAME{itbpF}{0.1859in}{0.1565in}{0.0138in}{}{}{t211-21-2.eps}{\raisebox{-0.1565in}{\includegraphics[height=0.1565in]{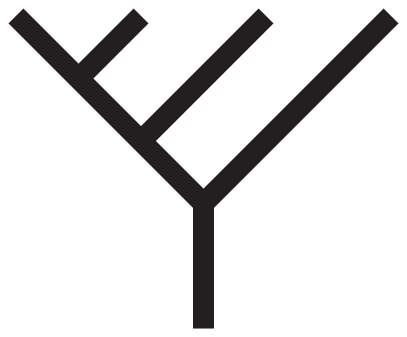}}}
$,$
\FRAME{itbpF}{0.1859in}{0.1565in}{0.0069in}{}{}{t121-21-2.eps}{\raisebox{-0.1565in}{\includegraphics[height=0.1565in]{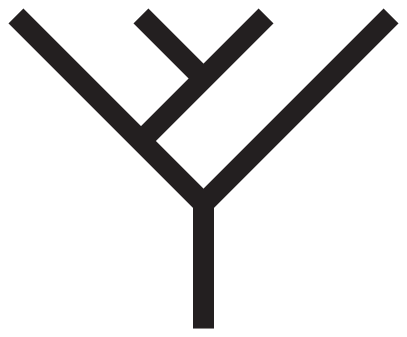}}}$\
,$
\FRAME{itbpF}{0.1859in}{0.1565in}{0.0138in}{}{}{t121-12-2.eps}{\raisebox{-0.1565in}{\includegraphics[height=0.1565in]{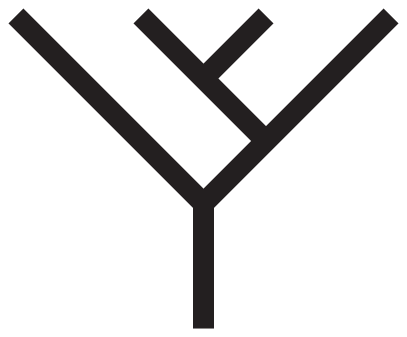}}} $,$ \FRAME{itbpF}{
0.1859in}{0.1565in}{0.0138in}{}{}{t22-2.eps}{\raisebox{-0.1565in}{\includegraphics[height=0.1565in]{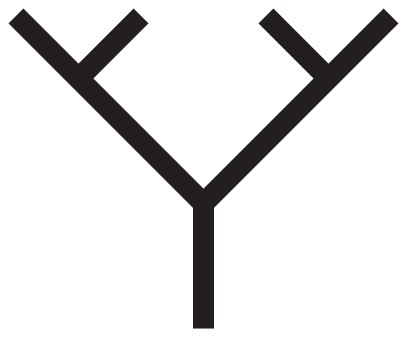}}}
and \FRAME{itbpF}{0.1859in}{0.1565in}{0.0138in}{}{}{t112-12-2.eps}{\raisebox{-0.1565in}{\includegraphics[height=0.1565in]{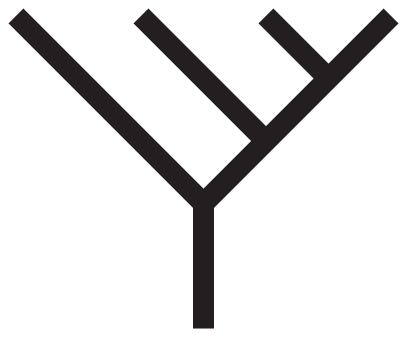}}} are DGA maps; \\
\FRAME{itbpF}{0.1859in}{0.1565in}{0.0277in}{}{}{t31-2.eps}{\raisebox{-0.1565in}{\includegraphics[height=0.1565in]{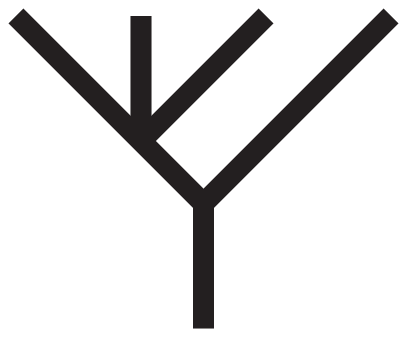}}} \ is a $($ \FRAME{itbpF}{%
0.1859in}{0.1565in}{0.0138in}{}{}{t211-21-2.eps}{\raisebox{-0.1565in}{\includegraphics[height=0.1565in]{T211-21-2.eps}}} $,$ \FRAME{itbpF}{%
0.1859in}{0.1565in}{0.0069in}{}{}{t121-21-2.eps}{\raisebox{-0.1565in}{\includegraphics[height=0.1565in]{T121-21-2.eps}}} $)$-derivation homotopy;
\\
\FRAME{itbpF}{0.1859in}{0.1565in}{0.0138in}{}{}{t121-3.eps}{\raisebox{-0.1565in}{\includegraphics[height=0.1565in]{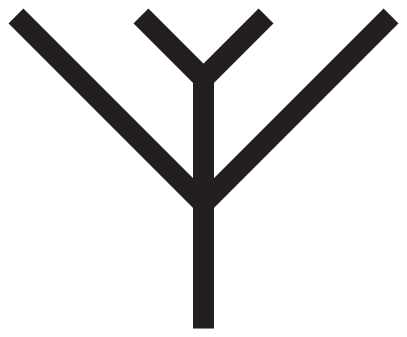}}} \ is a $($ \FRAME{itbpF}{%
0.1859in}{0.1565in}{0.0069in}{}{}{t121-21-2.eps}{\raisebox{-0.1565in}{\includegraphics[height=0.1565in]{T121-21-2.eps}}} $,$ \FRAME{itbpF}{%
0.1859in}{0.1565in}{0.0138in}{}{}{t121-12-2.eps}{\raisebox{-0.1565in}{\includegraphics[height=0.1565in]{T121-12-2.eps}}} $)$-derivation homotopy;
\\
\FRAME{itbpF}{0.1859in}{0.1565in}{0.0104in}{}{}{t13-2.eps}{\raisebox{-0.1565in}{\includegraphics[height=0.1565in]{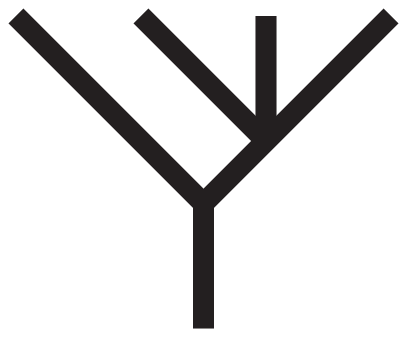}}} \ is a $($ \FRAME{itbpF}{%
0.1859in}{0.1565in}{0.0138in}{}{}{t121-12-2.eps}{\raisebox{-0.1565in}{\includegraphics[height=0.1565in]{T121-12-2.eps}}} $,$ \FRAME{itbpF}{%
0.1859in}{0.1565in}{0.0138in}{}{}{t112-12-2.eps}{\raisebox{-0.1565in}{\includegraphics[height=0.1565in]{T112-12-2.eps}}} $)$-derivation homotopy;
\\
\FRAME{itbpF}{0.1859in}{0.1565in}{0.0104in}{}{}{t211-3.eps} {\raisebox{-0.1565in}{\includegraphics[height=0.1565in]{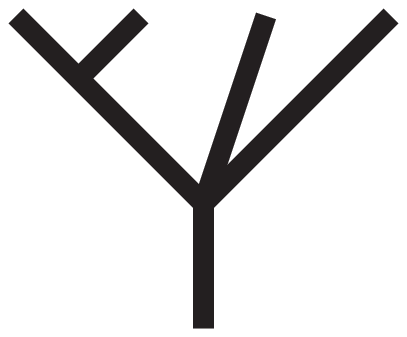}}} \ is a $($ \FRAME{itbpF}{%
0.1859in}{0.1565in}{0.0138in}{}{}{t211-21-2.eps}{\raisebox{-0.1565in}{\includegraphics[height=0.1565in]{T211-21-2.eps}}} $,$ \FRAME{itbpF}{%
0.1859in}{0.1565in}{0.0138in}{}{}{t22-2.eps}{\raisebox{-0.1565in}{\includegraphics[height=0.1565in]{T22-2.eps}}}
$)$-derivation homotopy; \\
\FRAME{itbpF}{0.1859in}{0.1565in}{0.0156in}{}{}{t112-3.eps}{\raisebox{-0.1565in}{\includegraphics[height=0.1565in]{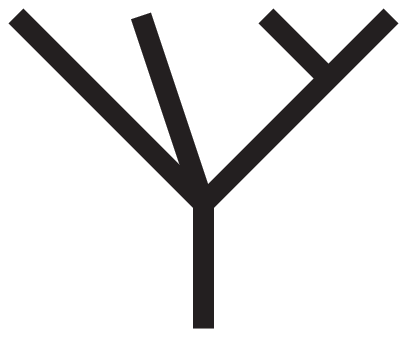}}} \ is a $($ \FRAME{itbpF}{%
0.1859in}{0.1565in}{0.0138in}{}{}{t22-2.eps}{\raisebox{-0.1565in}{\includegraphics[height=0.1565in]{T22-2.eps}}}
$,$ \FRAME{itbpF}{0.1859in}{0.1565in}{0.0138in}{}{}{t112-12-2.eps}{\raisebox{-0.1565in}{\includegraphics[height=0.1565in]{T112-12-2.eps}}} $)$-derivation homotopy.%
\end{tabular}%
\ \ \ \ \ \ \ \
\end{equation*}%
If \FRAME{itbpF}{0.1859in}{0.1565in}{0.0208in}{}{}{t4.eps}{\raisebox{-0.1565in}{\includegraphics[height=0.1565in]{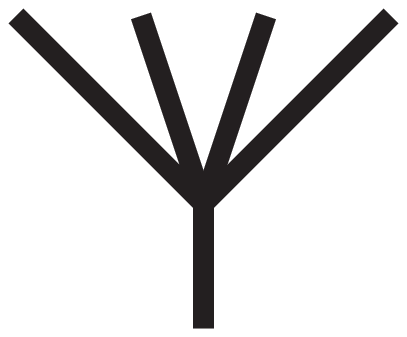}}} is a $\Delta $-derivation
homotopy with respect to $\mathfrak{F}_{4}$, then
\begin{equation*}
d^{\otimes 4}\text{ \FRAME{itbpF}{0.1859in}{0.1565in}{0.0208in}{}{}{t4.eps}{%
\raisebox{-0.1565in}{\includegraphics[height=0.1565in]{T4.eps}}} }-\text{ }\FRAME{itbpF}{%
0.1859in}{0.1565in}{0.0208in}{}{}{t4.eps}{\raisebox{-0.1565in}{\includegraphics[height=0.1565in]{T4.eps}}}%
\text{ }d=\left( \text{ \ \FRAME{itbpF}{0.1859in}{0.1565in}{0.0277in}{}{}{%
t31-2.eps}{\raisebox{-0.1565in}{\includegraphics[height=0.1565in]{T31-2.eps}}}
\ }+\text{ \ \FRAME{itbpF}{0.1859in}{0.1565in}{0.0138in}{}{}{t121-3.eps}{%
\raisebox{-0.1565in}{\includegraphics[height=0.1565in]{T121-3.eps}}} \ }+\text{ \ \FRAME{itbpF}{%
0.1859in}{0.1565in}{0.0104in}{}{}{t13-2.eps}{\raisebox{-0.1565in}{\includegraphics[height=0.1565in]{T13-2.eps}}}
\ }\right) -\left( \text{ \ \FRAME{itbpF}{0.1859in}{0.1565in}{0.0156in}{}{}{%
t112-3.eps}{\raisebox{-0.1565in}{\includegraphics[height=0.1565in]{T112-3.eps}}}
\ }+\text{ \ \FRAME{itbpF}{0.1859in}{0.1565in}{0.0104in}{}{}{t211-3.eps}{%
\raisebox{-0.1565in}{\includegraphics[height=0.1565in]{T211-3.eps}}} \ }\right) \text{ \ and}
\end{equation*}%
\begin{equation*}
\begin{tabular}{l}
\FRAME{itbpF}{0.1859in}{0.1565in}{0.0208in}{}{}{t4.eps}{\raisebox{-0.1565in}{\includegraphics[height=0.1565in]{T4.eps}}} $\mu =\mu ^{\otimes 4}$ $%
\sigma _{4,2}\left( \text{ \FRAME{itbpF}{0.1859in}{0.1565in}{0.0138in}{}{}{%
t211-21-2.eps}{\raisebox{-0.1565in}{\includegraphics[height=0.1565in]{T211-21-2.eps}}}}\ \otimes \text{ }\FRAME{itbpF}{0.1859in}{0.1565in}{0.0208in}{}{}{%
t4.eps}{\raisebox{-0.1565in}{\includegraphics[height=0.1565in]{T4.eps}}}%
\text{ }+\text{ }\FRAME{itbpF}{0.1859in}{0.1565in}{0.0208in}{}{}{t4.eps}{%
\raisebox{-0.1565in}{\includegraphics[height=0.1565in]{T4.eps}}}\text{ }\otimes \text{ \FRAME{%
itbpF}{0.1859in}{0.1565in}{0.0138in}{}{}{t112-12-2.eps}{\raisebox{-0.1565in}{\includegraphics[height=0.1565in]{T112-12-2.eps}}} \ }+\text{ \ }\FRAME{%
itbpF}{0.1859in}{0.1565in}{0.0242in}{}{}{t121-3.eps}{\raisebox{-0.1565in}{\includegraphics[height=0.1565in]{T121-3.eps}}}\text{ }\otimes \text{
\FRAME{itbpF}{0.1859in}{0.1565in}{0.0173in}{}{}{t13-2.eps}{\raisebox{-0.1565in}{\includegraphics[height=0.1565in]{T13-2.eps}}} \ }+\text{ \ \FRAME{itbpF}{%
0.1859in}{0.1565in}{0.0173in}{}{}{t211-3.eps}{\raisebox{-0.1565in}{\includegraphics[height=0.1565in]{T211-3.eps}}} }\otimes \text{ \FRAME{%
itbpF}{0.1859in}{0.1565in}{0.0208in}{}{}{t112-3.eps}{\raisebox{-0.1565in}{\includegraphics[height=0.1565in]{T112-3.eps}}}}\right. $ \\
\\
$\hspace{2.1in}\left. +\text{ \ }\FRAME{itbpF}{0.1859in}{0.1565in}{0.0277in}{%
}{}{t31-2.eps}{\raisebox{-0.1565in}{\includegraphics[height=0.1565in]{T31-2.eps}}}%
\text{ }\otimes \text{ }\FRAME{itbpF}{0.1859in}{0.1565in}{0.0208in}{}{}{%
t121-3.eps}{\raisebox{-0.1565in}{\includegraphics[height=0.1565in]{T121-3.eps}}}%
\text{ \ }-\text{ \ }\FRAME{itbpF}{0.1859in}{0.1565in}{0.0277in}{}{}{%
t31-2.eps}{\raisebox{-0.1565in}{\includegraphics[height=0.1565in]{T31-2.eps}}}%
\text{ }\otimes \text{ \FRAME{itbpF}{0.1859in}{0.1565in}{0.0173in}{}{}{%
t13-2.eps}{\raisebox{-0.1565in}{\includegraphics[height=0.1565in]{T13-2.eps}}}
}\right) ,$%
\end{tabular}%
\ \ \ \ \ \ \ \
\end{equation*}
and so on.

Note that when the initial family of maps is $\{\Delta=\Delta_2,\Delta_n\}$, the only components of $\Delta_{K}$ that come into play are primitive and the $\Delta_{K}$-compatible family $\mathfrak{F}_n$ is simply 
$$\{ f= ( \Delta \otimes \mathbf{1}^{\otimes n-2})\cdots (\Delta \otimes \mathbf{1})\Delta,\ \ g=( \mathbf{1}^{\otimes n-2}\otimes \Delta )\cdots (\mathbf{1} \otimes \Delta)\Delta\}.$$ We shall refer to a $\Delta$-derivation with respect to such an $\mathfrak{F}_n$ as an ``$(f,g)$-derivation of degree $n-2$.''
\end{example}

\begin{definition}
Let $(A,\mu_A)$ and $(B,\mu_B)$ be graded algebras and let $f,g:A\rightarrow B$ be algebra maps. A map $h:A\rightarrow B$ of degree $k$ is an $(f,g)$\textbf{-derivation of degree $k$} if 
$$h\mu_A=\mu_B(f\otimes h+h\otimes g).$$
\end{definition}
    
\section{The Hopf $A_{\infty}$-coalgebra $E\otimes\Gamma$}

The natural $A_\infty$-bialgebra structure on tensor factors $E \otimes \Gamma$ of $H_*(\mathbb{Z},n;\mathbb{Z}_p)$ is what we shall call a \textquotedblleft Hopf $A_\infty$-coalgebra.\textquotedblright \ Let $$f^{n}=\left(\Delta\otimes\mathbf{1}^{\otimes n-2}\right)  \cdots\left(  \Delta
\otimes\mathbf{1}\right)  \Delta.$$

\begin{definition}
\label{Hopf}Let $n\geq3.$ A \textbf{Hopf }$A_{\infty}$\textbf{-coalgebra} is an
$R$-Hopf algebra $(  A,\Delta,\mu)  $ together with an operation
$\Delta_{n}\in Hom^{n-2}(  A,A^{\otimes n})  $ such that 
$(A,\Delta,\Delta_{n})  $ is an $A_{\infty}$-coalgebra and $\Delta_{n}$ is
an $(  f^{n},f^{n})  $-derivation of degree $n-2$.
\end{definition}

\noindent Our choice of terminology in Definition \ref{Hopf} is motivated by
the fact that the operations $\left\{  \Delta,\Delta_{n},\mu\right\}  $
satisfy the Hopf relation and its analogue in degree $n-2$, i.e,
\begin{equation}\label{compatible}
\Delta_{n}\mu=\mu^{\otimes n}\sigma_{n,2}(f^{n}\otimes\Delta_{n} + \Delta_{n}\otimes f^{n}).
\end{equation}

\bigskip

Our main result applies the following lemma, which follows from Vandermonde's Identity:

\begin{theorem}
[Vandermonde's Identity]For $r,s\geq0$ and $0\leq k\leq r+s$,
\[
\binom{r+s}{k}=\sum_{i=0}^{k}\binom{r}{i}\binom{s}{k-i}.
\]

\end{theorem}

\begin{lemma}
\label{comb} Let $R=\mathbb{N\cup}\left\{  0\right\}  $ or $R=\mathbb{Z}_{p}$
with $p$ prime. For all $i\geq0$ and all $n$-tuples $(z_{1},\dots,z_{n})\in
R^{n}$ we have
\[
\binom{z_{1}+\cdots+z_{n}+1}{i}=\sum_{s_{1}+\cdots+s_{n}=i-1}\binom{z_{1}%
}{s_{1}}\cdots\binom{z_{n}}{s_{n}}+\sum_{t_{1}+\cdots+t_{n}=i}\binom{z_{1}%
}{t_{1}}\cdots\binom{z_{n}}{t_{n}},
\]
where we reduce $\operatorname{mod}p$ when $R=\mathbb{Z}_{p}.$ $\bigskip$
\end{lemma}

\begin{proof}
A standard formula for binomial coefficients gives
\[
\binom{z_{1}+\cdots+z_{n}+1}{i}=\binom{z_{1}+\cdots+z_{n}}{i-1}+\binom
{z_{1}+\cdots+z_{n}}{i}.
\]
Iteratively apply the Vandermonde's identity to the right-hand summand and
obtain
\begin{align}
\binom{z_{1}+\cdots+z_{n}}{i} =\sum_{k_{1}=0}^{i}\sum_{k_{2}=0}^{k_{1}}%
\hspace*{-0.05in}\cdots\hspace*{-0.05in}\sum_{k_{n-1}=0}^{k_{n-2}}\binom
{z_{1}}{k_{n-1}}\binom{z_{2} }{k_{n-1}-k_{n-2}}\cdots\binom{z_{n}}{i-k_{1}}.
\label{binomial}%
\end{align}
Note that the sum of the lower entries in the $n$ binomial coefficients of this
last expression is $i\ $and set $$t_{1}=k_{n-1},\ t_{2}=k_{n-1}-k_{n-2},
\cdots,\ t_{n-1}=k_{1}-k_{2}\textrm{ and }t_{n}=i-k_{1}.$$ Then expression
(\ref{binomial}) can be rewritten as
\[
\binom{z_{1}+\cdots+z_{n}}{i}=\sum_{t_{1}+\cdots+t_{n}=i}\binom{z_{1}}{t_{1}
}\cdots\binom{z_{n}}{t_{n}}.
\]

\end{proof}

\begin{remark}
The formula in Lemma \ref{comb} counts the number of ways $i$ objects can be
selected from a collection of $z_{1}+\cdots+z_{n}+1$ objects of $n+1$
different colors, one of which is uniquely colored \textquotedblleft
black\textquotedblright\ and $z_{i}$ of which have the same unique color for
each $i$. The first sum on the right-hand side counts the ways to select $i$
objects one of which is black; the second sum counts the ways to select $i$
objects none of which are black.
\end{remark}

\begin{theorem}
\label{EtimesG}Let $p$ be an odd prime. For each $i\geq0$ and $n\geq1,$ let
\[
A_{i}=E(v_{i},2np^{i}+1)\otimes\Gamma(w_{i},2np^{i+1}+2),
\]
let $\mu=(\mu_{E}\otimes\mu_{\Gamma})\sigma_{2,2},$ and for $j=2,p,$ let
$\Delta_{j}$ be defined as in (\ref{Delta2}) and (\ref{Deltap}). Then $\left(
A_{i},\Delta_{2},\Delta_{p},\mu\right)  $ is a Hopf $A_{\infty}$-coalgebra
over $\mathbb{Z}_{p}$.
\end{theorem}

\begin{proof}
For each $n\geq1$ and $i\geq0,$ $\left(  A_{i},\Delta_{2},\mu\right)  $ is a
Hopf algebra via $\mu=(\mu_{E}\otimes\mu_{\Gamma})\sigma_{2,2},$ where
\[
\mu_{E}(v_{i}\otimes v_{i})=0\text{ \ and }\mu_{\Gamma}(\gamma_{i}\left(
w_{i}\right)  \otimes\gamma_{j}\left(  w_{i}\right)  )=\binom{i+j}{i}%
\gamma_{i+j}\left(  w_{i}\right)  ,
\]
and $\left(  A_{i},\Delta_{2},\Delta_{p}\right)  $ is an $A_{\infty}%
$-coalgebra by Theorem \ref{berciano}. Furthermore, both sides of relation
(\ref{compatible}) vanish on tensor products involving $v_{i},$ and it is
sufficient to evaluate this relation on $\gamma_{i}\otimes\gamma_{j}%
=\gamma_{i}\left(  w_{i}\right)  \otimes\gamma_{j}\left(  w_{i}\right)  $.
First,
\begin{align}
\Delta_{p}\mu\left(  \gamma_{i}\otimes\gamma_{j}\right)   &  =\binom{i+j}%
{i}\Delta_{p}(\gamma_{i+j})\nonumber\\
&  =\sum_{z_{1}+\cdots+z_{p}=i+j-1}\binom{i+j}{i}v_{i}\gamma_{z_{1}}%
\otimes\cdots\otimes v_{i}\gamma_{z_{p}}\nonumber\\
&  =\sum_{z_{1}+\cdots+z_{p}=i+j-1}\binom{z_{1}+\cdots+z_{p}+1}{i}u,
\label{one}%
\end{align}
\newline where $u=v_{i}\gamma_{z_{1}}\otimes\cdots\otimes v_{i}\gamma_{z_{p}%
},$ and second,
\begin{align}
\mu^{\otimes p}  &  \sigma_{p,2}\left(  f^{p}\otimes\Delta_{p}+\Delta
_{p}\otimes f^{p}\right)  \left(  \gamma_{i}\otimes\gamma_{j}\right)
\nonumber\label{two}\\
&  =\sum_{\overset{l_{1}+\cdots+l_{p}=j}{s_{1}+\cdots+s_{p}=i-1}}\binom
{l_{1}+s_{1}}{s_{1}}\cdots\binom{l_{p}+s_{p}}{s_{p}}v_{i}\gamma_{l_{1}+s_{1}%
}\otimes\cdots\otimes v_{i}\gamma_{l_{p}+s_{p}}\nonumber\\
&  +\sum_{\overset{m_{1}+\cdots+m_{p}=j-1}{t_{1}+\cdots+t_{p}=i}}\binom
{m_{1}+t_{1}}{t_{1}}\cdots\binom{m_{p}+t_{p}}{t_{p}}v_{i}\gamma_{m_{1}+t_{1}%
}\otimes\cdots\otimes v_{i}\gamma_{m_{p}+t_{p}}\nonumber\\
&  =\sum_{z_{1}+\cdots+z_{p}=i+j-1}\left[  \sum_{s_{1}+\cdots+s_{p}=i-1}%
\binom{z_{1}}{s_{1}}\cdots\binom{z_{p}}{s_{p}}+\sum_{t_{1}+\cdots+t_{p}%
=i}\binom{z_{1}}{t_{1}}\cdots\binom{z_{p}}{t_{p}}\right]  u.
\end{align}
But expressions (\ref{one}) and (\ref{two}) are equal modulo $p$ by Lemma
\ref{comb}.
\end{proof}

In view of Eilenberg and Mac Lane's decomposition of $H_{\ast}\left(
\mathbb{Z},n;\mathbb{Z}_{p}\right)  $ discussed in Section 3, the submodule
$$A_{i}=E(v_{i},2np^{i}+1)\otimes\Gamma(w_{i},2np^{i+1}+2)\subset H_{\ast
}\left(  \mathbb{Z},n;\mathbb{Z}_{p}\right)  $$ is a naturally occurring Hopf
$A_{\infty}$-coalgebra for each $i\geq0$. Furthermore, when $n=3$, the Hopf
$A_{\infty}$-coalgebra structures on the $A_{i}$'s induce a global $A_{\infty
}$-bialgebra structure on $H_{\ast}\left(  \mathbb{Z},3;\mathbb{Z}_{p}\right)
$ as a formal (as yet undefined) tensor product of $A_{\infty}$-bialgebras.
The relationship between this global structure and the underlying topology of
$K\left(  \mathbb{Z},3\right)  $ is an interesting open question.

\end{document}